\documentclass[12pt]{article}

\usepackage{amsmath, amsfonts, amssymb}
\usepackage{graphicx}
\usepackage{float}
\usepackage{authblk}
\usepackage{caption}
\usepackage{subcaption}
\usepackage{url}
\usepackage{xcolor}
\definecolor{mygreen}{RGB}{0,128,0}

\usepackage{hyperref}

\title{From Obstacle Problems to Neural Insights: Feed Forward Neural Network Modeling of Ice Thickness}

\author[1]{Kapil Chawla\thanks{Email: kapchaw@iu.edu}}
\author[2]{William Holmes\thanks{Email: wrholmes@iu.edu}}
\author[1]{Roger Temam\thanks{Email: temam@indiana.edu}}
\affil[1]{Department of Mathematics and the Institute for Scientific Computing and Applied Mathematics, Indiana University, Indiana 47405, USA}
\affil[2]{Department of Mathematics and Cognitive Science Program, Indiana University, Indiana 47405, USA}
\affil[*]{Corresponding author: Kapil Chawla, kapchaw@iu.edu}

\begin{document}

\maketitle
\begin{abstract}
In this study, we integrate the established obstacle problem formulation from ice sheet modeling \cite{Ed2012, PT2023} with cutting-edge deep learning methodologies to enhance ice thickness predictions, specifically targeting the Greenland ice sheet. By harmonizing the mathematical structure with an energy minimization framework tailored for neural network approximations, our method's efficacy is confirmed through both 1D and 2D numerical simulations. Utilizing the NSIDC-0092 dataset for Greenland \cite{Bamber2001} and incorporating bedrock topography for model pre-training, we register notable advances in prediction accuracy. Our research underscores the potent combination of traditional mathematical models and advanced computational techniques in delivering precise ice thickness estimations.
\end{abstract}

% Keywords

% Introduction
\section{Introduction}
The melting of ice sheets, driven by climate change, is a topic of mounting concern across various scientific disciplines. This phenomenon is pivotal for understanding the dynamic processes of Earth's climate, particularly in regions such as Greenland. The melting of Greenland's ice sheet not only contributes to global sea level rise but also provides insights into intricate climate interactions and feedback loops. As such, mathematicians have developed complex models to delve deeper into ice sheet dynamics. Among these models, obstacle problems \cite{JF1987} offer a unique lens, presenting challenges in partial differential equations (PDEs). Over the years, numerous numerical methods have been devised to address these challenges. Most of these methods focus on providing approximation solutions to the weak variational inequality. Techniques like the Galerkin least squares finite element method (\cite{Xia2007}, \cite{Hug1992}, \cite{Erik2017}), multigrid algorithm (\cite{R1994}, \cite{R1996}), piecewise linear iterative algorithm \cite{Yuan2012}, the first-order least-squares method \cite{Tom2012}, the level set method \cite{k2004}, and the dynamical functional particle method \cite{R2020} have been employed with varying degrees of success.

In the wake of technological advancements, deep learning has emerged as a promising tool in many scientific applications. It has recently gained significant traction in solving differential equations and inverse problems (\cite{Rudy2019}, \cite{Raissi2019}, \cite{E2018}, \cite{Nick2019}, \cite{Rohit2018}, \cite{Khoo2019}). Despite this momentum, its application to variational inequalities remains in its infancy. Some studies (\cite{Cheng2023}, \cite{R1984}) have innovatively applied deep learning to the traditional obstacle problem, whereas others \cite{Huang2021} have ventured into using deep learning techniques for elliptic hemivariational inequalities. A common observation, however, is that many of these studies prioritize computations over theoretical insights.

With this backdrop, our paper endeavors to bridge this gap. We explore both traditional ice-sheet models \cite{Ed2012, PT2023} and introduce a computational approach using deep learning to address the obstacle problem, deriving inspiration from its variational form. A central theme of our work is to discern the influence of parameters such as network size and training samples on the outcomes. Through rigorous numerical experiments, we substantiate the efficacy of our proposed method.

This article is organized as follows: Section 2 introduces the mathematical formulation of the model and elucidates the ice-thickness variational inequality. Section 3 sheds light on the energy minimization formulation. Section 4 delineates the approximation of the solution using fully connected feedforward deep neural networks, detailing its architecture, universality as an approximator, and the composite loss function tailored for optimal training and optimization. Section 5 showcases numerical experiments for one and two-dimensional problems, accompanied by solution visualizations and error analysis. Section 6 applies our model to data sourced from Greenland. We wrap up in Section 7, offering a concise summary of our study's principal insights and findings.

% Methodology
\section{Mathematical Formulation of the Model}

In this section, we present the mathematical formulation of the model, which is adapted from the work presented in \cite{Ed2012, PT2023}. 

Let  \(\mathbb{R}^n\) denote the n-dimensional Euclidean space, equipped with the standard Euclidean norm. A domain \( \Omega \) in \( \mathbb{R}^n \) is defined as a bounded and connected open subset of \( \mathbb{R}^n \), whose boundary is Lipschitz continuous. Consider a subset \( \Omega \) residing within \( \mathbb{R}^2 \). For any point \( x = (x_1, x_2) \) contained within the closure of \( \Omega \), denoted as \( \bar{\Omega} \), we will utilize common mathematical operators without going into their detailed definitions here.

The bedrock elevation is denoted by the function \( b: \Bar{\Omega} \rightarrow \mathbb{R} \). It's noteworthy that a positive value of \( b \) represents elevations above sea level, while negative values correspond to depths below the sea level.

\begin{figure}[ht]
  \centering
  \includegraphics[width=0.7\textwidth]{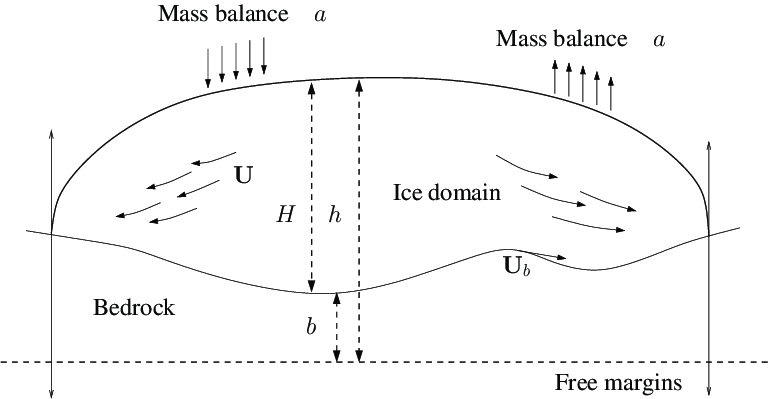}
  \caption{Cross-sectional view of an ice sheet with the respective notation, based on Jouvet et al.~(2012)}
  \label{fig:Ed_diag}
\end{figure}

Similarly, the elevation of the top surface of the ice sheet is characterized by the function \( h: \Bar{\Omega}\rightarrow \mathbb{R} \). It is imperative to emphasize that throughout the domain \( \Omega \), \( h \) always maintains a value greater than or equal to \( b \). A visual representation of this relationship is provided in Figure~\ref{fig:Ed_diag}. Consequently, the thickness of the ice, denoted as \(H: h-b\), consistently remains nonnegative throughout \( \Bar{\Omega} \). This insight underscores the observation that studying changes in ice thickness is tantamount to addressing an obstacle problem, where the bedrock acts as the primary constraint.

This particular constraint implies the existence of a free boundary. Let's define 
\begin{equation}
\Omega^{+} = \{h>b\}=\{H>0\}
\end{equation}
as the region within \( \Omega \) where the ice is present. Accordingly, the associated free boundary can be expressed as:
\begin{equation}
    \Gamma_f = \Omega \cap \partial \Omega^{+}.
\end{equation}

For a continuous ice sheet, wherein \(\Omega^+\) is populated, the source function \(a\) should inherently display a positive value, denoting ice accumulation, within specified regions of \(\Omega^+\). However, considering a defined free margin and the continuity of \(a\), \(a\) is expected to yield a negative value, indicative of ablation, outside \(\Omega^+\), especially within \(\Omega^- = \Omega \backslash (\Omega^+ \cup \Gamma)\). It is in this dynamic that the ice migrates from accumulation zones to areas dominated by ablation, resulting in the ice sheet tapering to an imperceptible thickness at its edges.

Moreover, the ice thickness is also influenced by the basal sliding velocity, \( U_b \). This velocity can be conceptualized as a specific vector field within \( \mathbb{R}^2 \), which is contingent solely on its lateral position. This can be formally described as:
\[ U_b : \Bar{\Omega} \rightarrow \mathbb{R}^2. \]
As documented in references such as \cite{R2009}, it is generally accepted that, when the ice base is in a frozen state, \( U_b \) inherently equals zero.

The horizontal ice flow velocity is represented by the vector field \( U: \bar{\Omega} \times \mathbb{R} \rightarrow \mathbb{R}^2 \). The viscosity of the ice is described in terms of the Glen power law with ice softness coefficient \( A(x,z) \) and exponent \( 2.8 \leq p \leq 5 \). The attainable values for \( p \) are suggested by laboratory experiments \cite{lab}. The vector field \( U \) can be computed using the surface elevation \( h \) and its gradient as discussed in \cite{R2009}.

\begin{equation} \label{vel_U}
    U(x,z) = -(2\rho g)^{p-1} \left[ \int_b^z A(s) (h-s)^{p-1} ds\right] |\nabla h|^{p-2} \nabla h + U_b
\end{equation}
Here, \( \rho \) denotes the ice's density and \( g \) symbolizes the gravitational acceleration.

The morphology of ice sheets is independent of the surface equation, which correlates the ice surface's motion with the ice velocity and the mass balance data \( a=a(x,z) \) \cite{jouvet_rappaz_bueler_blatter_2011}. Equivalently for this incompressible flow, in steady state a continuity equation applies to the flow where the volume flux is characterized by the vector field \( \textbf{q} \), defined as the integral of the ice flow velocity \( U \) concerning the vertical direction \cite{R2009}:
\begin{equation}\label{flux_q}
     \nabla \cdot \textbf{q} = a
\end{equation}
Given by:
\begin{equation}\label{ice_eqn}
    \textbf{q}=\int_b^h U(z) dz
\end{equation}
By inclusion of (\ref{vel_U}) and (\ref{ice_eqn}) into (\ref{flux_q}), the equation for surface elevation \( h \) can be formulated as:
\begin{equation}\label{ice_orig_eqn_1}
    - \nabla \cdot \left( -(2\rho g)^{p-1} \left[ \int_b^z A(s) (h-s)^{p-1} ds\right] |\nabla h|^{p-2} \nabla h - (h-b) U_b \right)=a
\end{equation}
For a scenario where \( A(x,z) = A_0 \) remains constant, Equation (\ref{ice_orig_eqn_1}) simplifies to:
\begin{equation}\label{ice_orig_eqn_2}
    - \nabla \cdot \left( K (h-b)^{p+1} |\nabla h|^{p-2} \nabla h - (h-b) U_b\right) = a
\end{equation}
with \( K \) being a positive constant. It's imperative to note that Equations (\ref{ice_orig_eqn_1}) and (\ref{ice_orig_eqn_2}) are only applicable within the domain containing ice.

To form a weak solution for problem (\ref{ice_orig_eqn_2}) over the entire domain \( \Omega \), one must posit the absence of a volumetric ice flow towards \( \Omega^{-} \):
\begin{equation}\label{nu_0}
    \textbf{q}\cdot \nu = 0~~~\text{and}~~~H=0~~\text{on}~~\Gamma_f
\end{equation}
Here, \( \nu \) represents the outward unit normal vector along \( \Gamma \). Taking (\ref{nu_0}) into account, it's intuitive to extrapolate the volume flux \( \textbf{q} \) as zero outside \( \Omega^+ \):
\begin{equation}
    \textbf{q}=0~~\text{in}~~\Omega^{-}
\end{equation}
Let's introduce \( v \) as a smooth test function that satisfies \( v \geq b \) throughout \( \Bar{\Omega} \) and multiply (\ref{ice_eqn}) by \( (v-h) \), subsequently integrating over \( \Omega \). Applying Gauss's theorem yields:
\begin{equation}
    -\int_{\Omega} \textbf{q} \cdot \nabla (v-h) = \int_{\Omega^{-}} (\nabla \cdot \textbf{q})(v-h) + \int_{\Omega^{+}} (\nabla \cdot \textbf{q})(v-h) -  \int_{\Gamma_f} [(v-h) \textbf{q} \cdot \nu]_{-}^+
\end{equation}
Here, the difference across \( \Gamma_f \) is represented by \( [~]_{-}^{+} \). Given the conditions \( \textbf{q} = 0 \) and \( a \leq 0 \) within \( \Omega^- \), it follows that \( \nabla \cdot \textbf{q} \geq a \) in \( \Omega^- \). Moreover, since \( \textbf{q} \cdot \nu = 0 \) on \( \Gamma_f \), we deduce:
\begin{equation}
    - \int_\Omega \textbf{q} \cdot \nabla(v-h) \geq \int_\Omega a(v-h)
\end{equation}
Expressing \( \textbf{q} \) based on (\ref{ice_orig_eqn_2}), we encounter the subsequent variational inequality:
\begin{equation}\label{ineq}
    \left((2\rho g)^{p-1} \left[ \int_b^z A(s) (h-s)^{p-1} ds\right] |\nabla h|^{p-2} \nabla h - (h-b) U_b\right) \cdot \nabla (v-h) \geq \int_\Omega a(v-h)
\end{equation}
An essential observation to make is that \( \textbf{q} \) displays a degenerate behavior at the free boundary. To navigate the gradient degeneracy proximate to the free boundary, a transformation was introduced in \cite{calvo2002}:
\begin{equation}
    H = u^{(p-1)/(2p)}
\end{equation}
This results in a transformed version of Equation (\ref{ice_orig_eqn_2}):
\begin{equation}
    - \nabla \cdot \left( \mu(x,u) |\nabla u - \Phi(x,u)|^{p-2} \nabla u - (h-b) U_b\right) = a
\end{equation}
Where \( \mu(x,u) \) and \( \Phi(x,u) \) are specifically defined functions. Given this transformation, further study would reveal the implications and applications of the aforementioned equations.

In the present work, we focus on simplifying the nonlinear dynamics by making specific assumptions to address a more tractable problem. Specifically, we consider the scenario where \( U_b = 0 \). Under this assumption, we set \( \Psi(x,u) = 0 \) and treat \( \mu(x,u) \) as a constant, denoted \( \mu_0 \). By adopting these assumptions, the inequality expressed in (\ref{ineq}) simplifies to the following:

\begin{equation}\label{ineq_final}
    \int_\Omega \left( \mu_0 |\nabla u - \Phi(x,u)|^{p-2}(\nabla u - \Phi(x,u))\right) \cdot \nabla(v-u) \geq \int_\Omega a(v-u)
\end{equation}

This inequality holds for all \( v \geq 0 \) defined over \( \Omega \) and \( p \geq 2 \). For simplification, we take $\mu_0=1$.\\

In the seminal work by the original authors, several inherent nonlinearities were identified but not thoroughly explored. Recognizing the significance and implications of these nonlinearities, we have taken steps in this paper to address and analyze them more holistically. Specifically, we have chosen to focus on two primary nonlinearities that were presented but not fully addressed in the original study \cite{Ed2012}:

\begin{enumerate}
    \item A \(p\)-Laplace-type nonlinearity arising from Glen’s flow law. Notably, while the original authors noted the simplification of this nonlinearity when \(p = 2\), resulting in \(|\nabla u - \Phi(x, u)|^{p-2} = 1\), our investigation delves deeper into cases where \(p > 2\), thereby introducing and emphasizing the complexities of this nonlinearity.
    
    \item An intricate nonlinearity in \(\Phi(x, u)\) due to the bedrock gradient, \(\nabla b\). In our research, we consider scenarios where the bedrock elevation isn't flat, which introduces another layer of nonlinearity. This is in contrast to the original paper where the implications of a non-flat bedrock weren't fully explored, resulting in \(\Phi(x, u) = 0\).
\end{enumerate}

By focusing on these areas, our aim is to provide a more comprehensive understanding of the problem, shedding light on aspects that were previously left in the shadows. As we delve deeper into our discourse, our approach becomes particularly evident. We will explicitly address these nonlinearities, offering insights and solutions that not only build upon but also augment and enhance the foundation set by the original author, especially as we transition into the energy minimization formulation.

\section{Energy Minimization Formulation}

Consider a Hilbert space \( V \) such that \( V \in H^1(\Omega) \). Let the subspace \( U \) be defined by
\[
U = \left\{ v \in V : v \big|_{\partial \Omega} = 0 \right\}.
\]
Furthermore, we introduce the constraint set \( K \) as
\[
K = \left\{ v \in V : v \geq b~\text{in}~\Omega \text{ and } v \big|_{\partial \Omega} = 0 \right\}.
\]

Given these definitions, our variational inequality problem, as specified in (\ref{ineq_final}), can be equivalently posed as an energy minimization problem:
\begin{equation}\label{min}
\text{Find } u \in K \text{ such that } J(u) \leq J(v)~\text{for all } v \in K,
\end{equation}
with the energy functional \( J \) given by
\[
J(v) = B(v,v) - \langle a,v \rangle,
\]
where 
\[
B(v,v) = \frac{1}{p} \int_\Omega |\nabla v - \Psi|^p \, dx.
\]
We introduce a regularization term to the energy functional and define the augmented energy functional as
\begin{equation}\label{energy_funct}
    L(v) = J(v) + \alpha \int_\Omega [b(x) - v(x)]_+^2 \, dx,
\end{equation}
where \( \alpha \) is a positive constant and 
\[ [b(x) - v(x)]_+ = \max\{0, b(x) - v(x)\}. \]
Based on the findings from \cite{unique}, the unique solution \( u \) of problem (\ref{energy_funct}) approaches the solution of (\ref{min}) as \( \alpha \) grows significantly large. This leads us to the following minimization problem:
\begin{equation}
    \min_{v \in V} L(v),
\end{equation}
with the energy functional given by
\begin{equation}
    L(v) = \frac{1}{p} \int_\Omega |\nabla v - \Psi|^p \, dx - \int_\Omega a v \, dx + \alpha \int_\Omega [b - v]_+^2 \, dx,
\end{equation}
where \( \alpha \) is again a positive constant.

In numerical computations, this continuous energy functional is approximated by its discrete counterpart:
\begin{equation}
    \Bar{L}(v) = \frac{|\Omega|}{N} \sum_{i=1}^N \left[ \frac{1}{p} \|\nabla v(X_i) - \Psi(X_i)\|_2^p - a(X_i) v(X_i) + \alpha [b(X_i) - v(X_i)]_+^2 \right],
\end{equation}
where \( \{X_k\}_{k=1}^N \) are i.i.d random variables sampled from the uniform distribution \( U(\Omega) \).

If the solution on the boundary is non-zero, the minimization problem can be generalized to:
\begin{equation}
    \min_{v \in V} L(v),
\end{equation}
with the modified energy functional:
\begin{equation}
    L(v) = \frac{1}{p} \int_\Omega |\nabla v - \Psi|^p \, dx - \int_\Omega a v \, dx + \alpha \int_\Omega [b - v]_+^2 \, dx + \beta \int_{\partial\Omega} (v - h)^2,
\end{equation}
where \( h \) denotes the solution on the boundary.

With the energy functional clearly defined, the next step is to seek efficient computational or approximation methods for solutions. Here, Deep Neural Networks (DNNs) emerge as an especially promising tool. Their ability to capture intricate nonlinear relationships and represent functions make them apt for approximating our energy functional. In the following section, we delve into how we can leverage these networks for our purpose.

\section{Solution Approximation Using Deep Neural Networks}

We utilize a fully connected feedforward neural network, denoted as $\mathbf{f}: \mathbb{R}^d \rightarrow \mathbb{R}^{N_D}$, to approximate our solution. This neural network comprises multiple layers, with each layer introducing nonlinearity through an activation function. This ensures that the network can function as a universal approximator, capturing the complex relationships inherent to this problem. The weights and biases in the network are iteratively refined using backpropagation based on the minimization of a loss function associated with the energy minimization formulation above.

To measure the accuracy and efficacy of our deep neural network approximation, we employ a composite loss function derived from three primary components:

\begin{itemize}
    \item \textbf{Residual Loss (${loss_1}$)}: 
    \begin{equation}\label{loss_1}
    {loss_1} = \frac{1}{N} \sum_{i=1}^N \left[ \frac{1}{p}{\|\nabla u^*(X_i)-\Psi(X_i)\|^p} - a(X_i)u^*(X_i)\right].
    \end{equation}
    This term captures the degree to which the current state of the solution fails to satisfy the PDE. Minimizing this is equivalent to convergence to a solution of the PDE.

    \item \textbf{Obstacle Loss (${loss_2}$)}:
    \begin{equation}\label{loss_2}
    {loss_2} = \frac{1}{N} \sum_{i=1}^N \left[b(X_i) - u^{*}(X_i)\right]_{+}^2.
    \end{equation}
    This term penalizes the discrepancy between the network's output and the obstacle function, $b$. This component helps ensure the resulting solution lies above the constraint.

    \item \textbf{Boundary Condition Loss (${loss_3}$)}:
    \begin{equation}\label{loss_3}
    {loss_3} = \frac{1}{M} \sum_{j=1}^M \left[u^{*}(Y_j) - h(Y_j)\right]^2.
    \end{equation}
    This component ensures our solution adheres to known values on the boundary of our domain, effectively representing the boundary condition loss.
\end{itemize}

The total loss function is then a weighted combination of the aforementioned losses:

\begin{equation}\label{total_loss}
L(\theta) = {loss_1} + \alpha \cdot {loss_2} + \beta \cdot {loss_3}
\end{equation}

\noindent The weighting coefficients \( \alpha \) and \( \beta \) balance the importance of loss components for PDE satisfaction, constraint satisfaction, and boundary condition satisfaction respectively. Their optimal values, attained through iterative trial and error, ensure a balanced representation of all problem constraints in the loss function. Our training uses a batch approach, where data subsets iteratively refine the network parameters. With this foundational understanding, we now explore the neural network's practical applications in subsequent sections, showcasing its versatility and efficacy.

\section{Numerical Illustrations: 1D and 2D Examples}

Before describing the specific case studies, we outline the sequence of numerical experiments we will perform. Initially, we will use the "method of manufactured solutions" (MMS) to test this approach. This technique, commonly used in computational science, involves creating an exact solution (the "manufactured solution") for an augmented form of partial differential equation. This is done by choosing a solution function (that satisfies obstacle and boundary constraints), computing the residual arising from that function not satisfying the PDE, and augmenting that PDE with an inhomogeneity that offsets that residual exactly. The manufactured solution is then a solution to the inhomogeneity augmented PDE. The MMS allows for a comprehensive assessment of the numerical method in question as it provides a "ground truth" against which numerical approximations can be rigorously compared. By applying MMS in our preliminary experiments, we aim to validate the robustness and accuracy of this deep learning approach. 

This problem is comprised of both a complex PDE and an obstacle constraint. We will thus first apply this approach to a simpler PDE with an obstacle constratint followed by the more complex PDE of interest with the constraint. This will be performed on both 1D and 2D domains. Following this we will test this approach on the motivating problem of trying to find the elevation profile of ice above the bedrock of Greenland.

For the solution approximator, we adopt an architecture comprising a minimum of 5 hidden layers, each equipped with 128 neurons. We utilized the squared Rectified Linear Unit ($ReLU^2$) in our study to leverage its improved differentiability and empirical performance, particularly enhancing the accuracy and convergence in solving complex Partial Differential Equations (PDEs).Furthermore, to foster a consistent learning environment, layer normalization techniques are integrated into the model. The optimization phase employs the Adam variant of the stochastic gradient descent (SGD) method. The learning rate is initially set to $5 \times 10^{-4}$ for the first 500 iterations. After which, it reduces by half from the 500 to 750 iterations, and then further reduces by half for the remaining iterations. This learning rate schedule was determined through trial and error. The primary objective remains the minimization of the associated loss function. For the implementation and training processes, we relied on the capabilities offered by the Pytorch library \cite{pytorch}. We note that there are many variations on this network and training structure that could be investigated. Here we are mainly investigating whether this approach is feasible and refinement is left for future study.

\subsection{1D Example: Case with \texorpdfstring{$p=2$}{p=2}}

To initiate our numerical exploration, we consider a one-dimensional problem with $p=2$. For this instance, we use for the following obstacle function:
\begin{equation}
b(x) = 
\begin{cases}
10 \sin(2 \pi x), & \text{for } 0 \leq x \leq 0.25, \\
5 \cos( \pi (4x-1))+5, & \text{for } 0.25 \leq x \leq 0.75, \\
10 \sin(2 \pi (1-x)), & \text{for } 0.75 \leq x \leq 1.0.
\end{cases}
\end{equation}

\noindent The exact solution is characterized as:
\begin{equation}
u_{\text{exact}}(x) = 
\begin{cases}
10 \sin(2 \pi x), & \text{for } 0 \leq x \leq 0.25, \\
10, & \text{for } 0.25 \leq x \leq 0.75, \\
10 \sin(2 \pi (1-x)), & \text{for } 0.75 \leq x \leq 1.0.
\end{cases}
\end{equation}

\begin{figure}[H]
  \centering
  \includegraphics[width=0.7\textwidth]{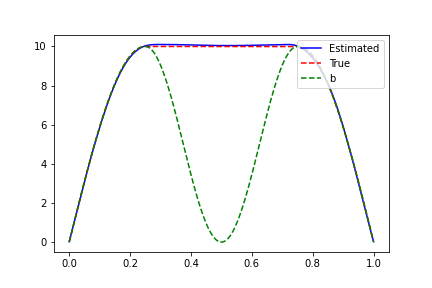}
  \caption{Comparison of the 1D solution, obstacle function (denoted \emph{b}), and the exact solution.}
  \label{fig:one_d_plot}
\end{figure}

\noindent Refer to Figure \ref{fig:one_d_plot} for the corresponding visualization.
\noindent To train a solution approximator for this problem, we follow the approach outlined above with appropriate modifications for the problems specifics. The points \( \{X_k\}_{k=1}^N \) for evaluation are drawn from a uniform distribution over the interval \( [0,1] \).

\begin{figure}[H]
    \centering
    \includegraphics[width= 0.9\textwidth]{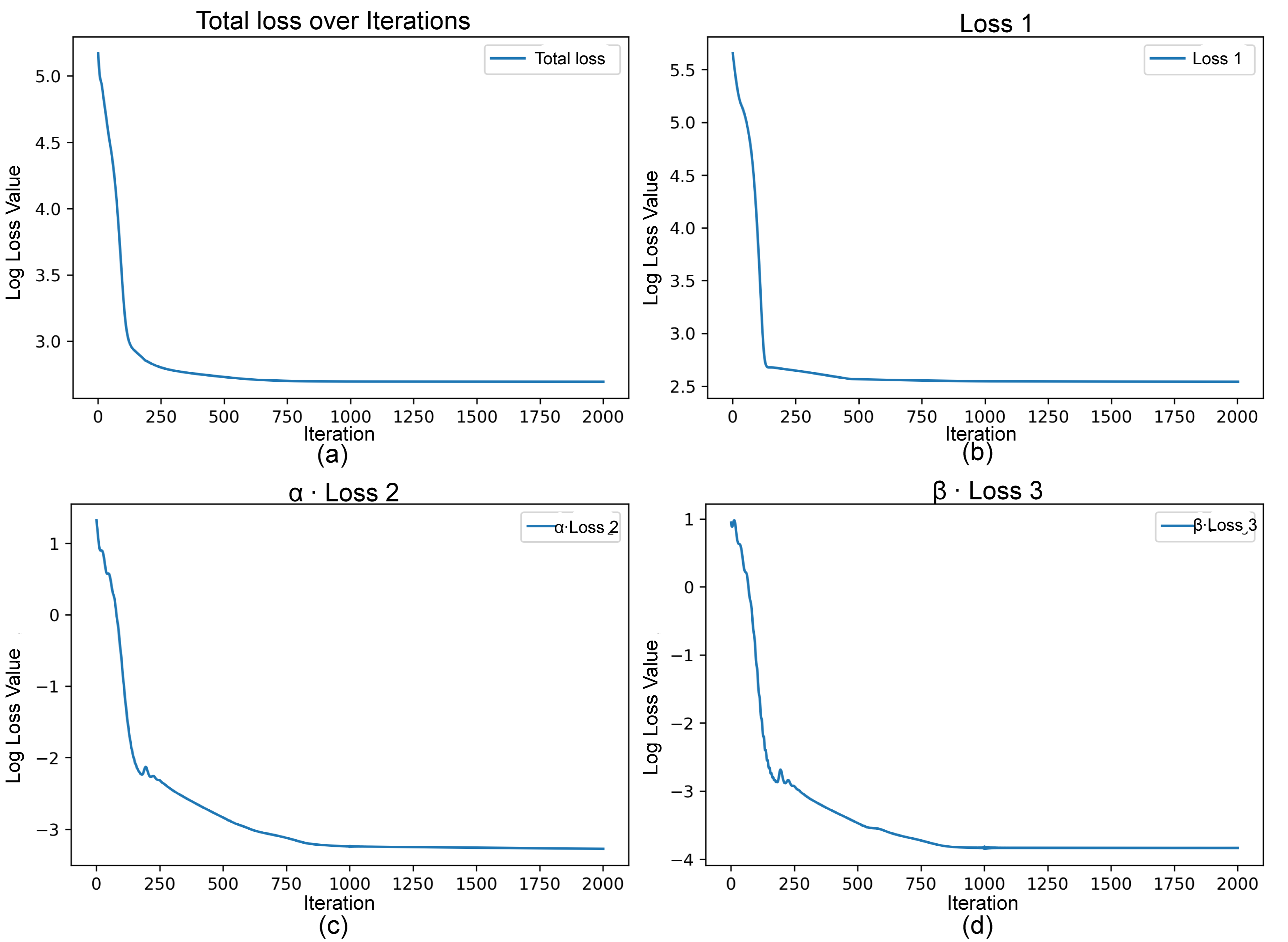}
    \caption{Evolution of different losses during training, plotted against the number of iterations on a logarithmic scale (with base 10). (a) \textit{Total Loss}, represented as ``Log Loss Value'' on the y-axis and computed from Equation \ref{total_loss}, symbolizes the amalgamation of all individual loss terms. (b) \textit{Loss 1}, denoted as ``Log Loss Value'' on the y-axis from Equation \ref{loss_1}, illustrates the residual loss term. (c) \textit{Loss 2}, weighted by coefficient $\alpha$ and articulated as ``Log Loss Value'' on the y-axis in Equation \ref{loss_2}, encapsulates the obstacle loss term. (d) \textit{Loss 3}, accentuated by coefficient $\beta$ and reflected as ``Log Loss Value'' on the y-axis from Equation \ref{loss_3}, highlights the boundary loss inherent in the problem.}

    \label{fig:lossplots_one_d}
\end{figure}

In our experimental setup, the neural network undergoes training over \( 2000 \) iterations with the regularization parameter \( \alpha \) and \(\beta\) fixed at \( 4000 \). We then assess the performance of our methodology on a mesh grid, uniformly spaced and marked by a fine resolution of \( 10^{-3} \). The resultant numerical solution is presented in Figure~\ref{fig:one_d_plot}, where it is juxtaposed against both the obstacle function and the exact solution. Delving deeper into the training dynamics, Figure~\ref{fig:lossplots_one_d} charts the evolution of various loss metrics, including the total loss, Loss 1, Loss 2 and Loss 3. A notable observation from our experiments is the efficiency of the proposed technique; it typically converges within a mere \( 500 \) iterations, underscoring its robustness and efficacy.

In our analysis, the difference between the exact and approximated solutions is quantified using the \(L^1\) norm, defined as:
\[ \text{loss}_{\text{exact}} = \frac{1}{N} \sum_{i=1}^{N} |U_{\text{inner\_pred},i} - U_{\text{exact\_inner},i}| \]
where \(N\) denotes the total number of grid points or samples, and \(U_{\text{inner\_pred},i}\) and \(U_{\text{exact\_inner},i}\) represent the approximated and exact solutions, respectively, at the \(i\)-th sample. A plot of this loss versus the iterations provides insight into the convergence behavior and accuracy of the neural network's predictions. These results demonstrate that this approach learns a neural network approximator that accurately learns to solve the PDE subject to the given constraints.

\begin{figure}[H]
  \centering
  \includegraphics[width=0.9\textwidth]{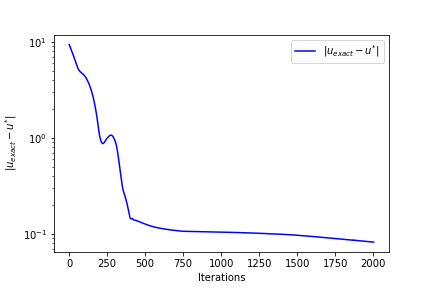}
  \caption{Evolution of the \(L^1\) error between the approximated solution \( u^* \) and the exact solution \( u_{\text{exact}} \) over training iterations. The plot showcases the difference in magnitudes as training progresses, shedding light on the convergence and accuracy of the neural network's predictions.}
  \label{fig:L_1_one_d}
\end{figure}
\subsubsection{Analysis of Relative Error in the Context of Sample Variability}

\begin{figure}[H]
    \centering
    \includegraphics[width=0.7\textwidth]{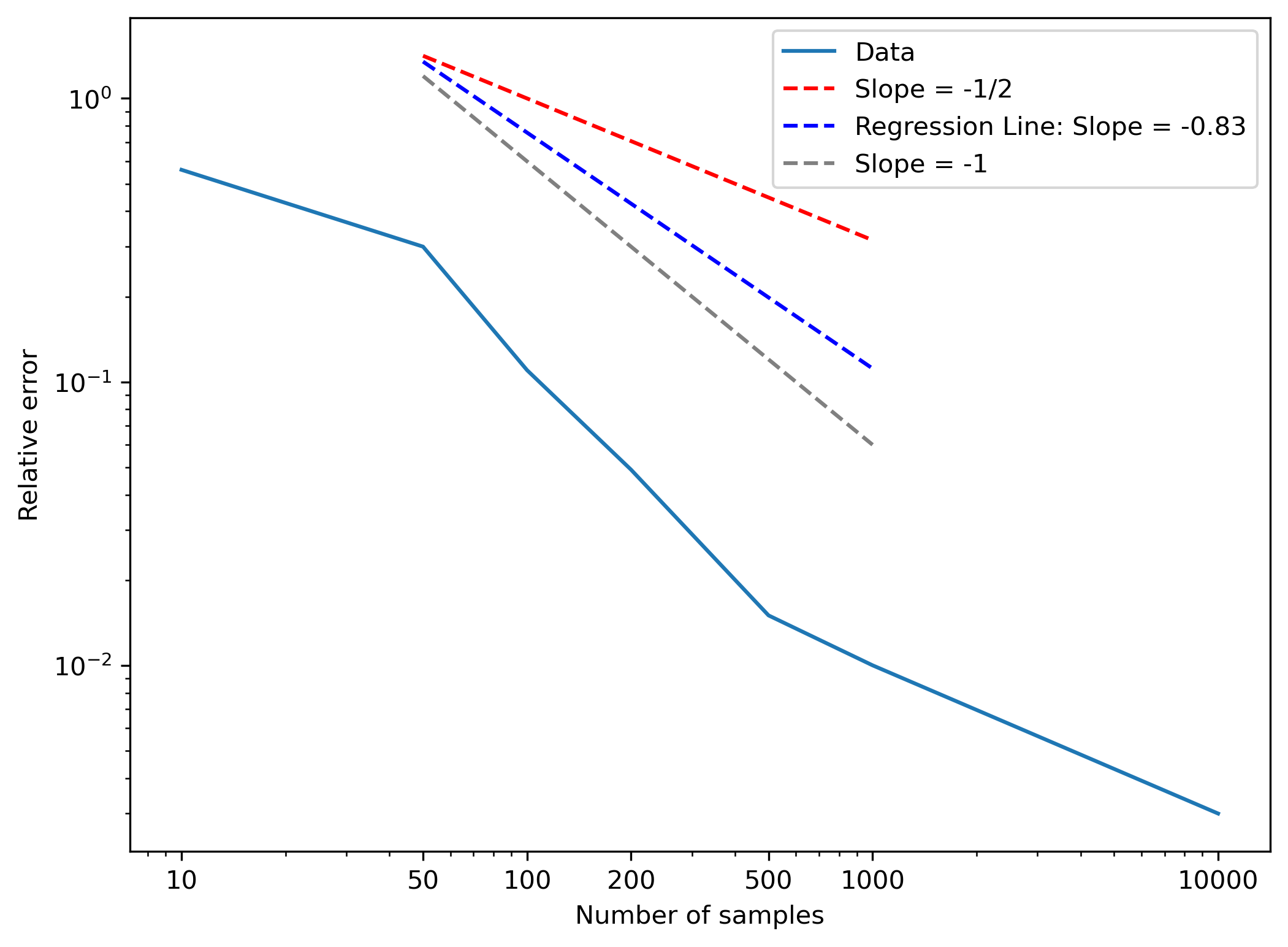}
    \caption{Graphical representation of relative error against varying sample sizes. The blue dashed line represents the regression line with a slope of -0.83, which lies between the theoretical slopes of -1/2 (red) and -1 (gray).}
    \label{fig:relative_error_plot}
\end{figure}

We next analyze the accuracy of this approach as a function of the number of sampled points at which the PDE loss is computed. This analysis provides insights into how the sampling granularity can influence the accuracy of the derived solution. These sample points are systematically chosen from a uniform grid, ensuring consistent and unbiased assessment. Figure~\ref{fig:relative_error_plot} portrays the relationship between sample size and the resultant relative error.
%{\color{blue}
In our analysis regarding the number of samples, Figure~\ref{fig:relative_error_plot} highlights the reduction in relative error as sample count augments. To better visualize the expected linear scaling of error with the number of points, a line with a -1 slope was added to the plot. This represents the ideal trend where, if the error scales as \( \frac{1}{N} \), then \( \log(\text{error}) \) directly corresponds to \( -\log(N) \).
%}
%{\color{red}In our analysis regarding the number of samples, Figure~\ref{fig:relative_error_plot} highlights the reduction in relative error as the sample count increases. To better visualize the expected scaling of error with the number of points, we included lines with slopes of -1 and -1/2 in the plot. The line with a -1 slope represents the ideal trend where the error scales as \( \frac{1}{N} \), corresponding to \( \log(\text{error}) \) directly relating to \( -\log(N) \cite{strikwerda2004finite}\). The line with a -1/2 slope represents the expected Monte Carlo error behavior, where the error scales as \( \frac{1}{\sqrt{N}} \cite{glasserman2004monte}\).

To further analyze the error scaling, we performed a linear regression on the log-transformed data points, resulting in a slope of approximately -0.83. This regression line, shown in the updated plot, lies between the theoretical slopes of -1 and -1/2. Therefore, the actual slope from our data indicates that the error reduction behavior in our analysis is intermediate between these two theoretical expectations.
%}

In the analysis (Figure~\ref{fig:relative_error_plot}) illustrates the reduction in relative error as the computational grid is refined.  Common computational PDE approaches (finite difference, finite element, etc.) produce errors that scale as 1/N \cite{strikwerda2004finite} while Monte Carlo approaches commonly lead to errors that scale as 1/sqrt(N) \cite{glasserman2004monte}. To visualize the scaling relationship between error and number of points (N) Figure 5 is plotted on log-log axes and a regression line is fit to the error magnitude to compute the scaling slope. Results indicate that errors decay faster than 1/sqrt(N) but more slowly than 1/N for the range used.

In our comprehensive analysis, we not only explored the impact of sampling but also assessed the influence of varying regularization parameters and loss weighting parameters \( \alpha \) and \( \beta \) on relative error. Understanding these dynamics is vital, as they shed light on the model's sensitivity to its central hyperparameters, ensuring that it captures underlying patterns effectively without the pitfalls of overfitting or underfitting.

Table~\ref{tab:relative_error_penalty} illustrates the relative errors for different values of the regularization parameter. These results illustrates that a proper balance between the different loss components is needed to ensure optimal training.

\begin{table}[H]
    \centering
    \begin{tabular}{|c|c|c|}
        \hline
        \textbf{Parameter $(\alpha)$} & \textbf{Parameter $(\beta)$} & \textbf{Relative Error} \\
        \hline
        100 & 100 &0.40 \\
        \hline
        500 & 100 &0.015 \\
        \hline
        1000 & 500 &0.0035 \\
        \hline
        4000 & 4000 &0.0010 \\
        \hline
        5000 & 4000&0.0023 \\
        \hline
        %... add as many rows as necessary
    \end{tabular}
    \caption{Relative errors corresponding to a spectrum of regularization parameters. }
    \label{tab:relative_error_penalty}
\end{table}

\subsection{2D Problem: General Case with Any \texorpdfstring{\( p \geq 2 \)}{p >= 2}}

We now consider the 2D problem with \( p \geq 2\). For this analysis, we use the collocation points that are uniformly sampled from the unit square. For simplicity and to facilitate the use of the MMS approach, we use a circularly symmetric obstacle function:
\begin{equation}
    b(r) = \begin{cases}
1 - \left( r^{\frac{p}{p-1}} - (1-r)^{\frac{p}{p-1}} + 1 - \frac{p}{p-1}r \right) & \text{for } r \leq r^*\\
0 & \text{for } r > r^* \\
\end{cases}
\label{eq:obstacle_function}
\end{equation}
where \( r \) represents the radial distance in our 2D space and \( r^* = 0.75 \) is a threshold value defining the boundary of the obstacle's influence.

We construct an exact solution within this domain \( \Omega = [0,1] \times [0,1] \) given by:

\begin{equation}
    u_{exact}(r) = 
    \begin{cases}
        1 - \left( F(r) - G(r) + 1 - E(r) \right) & \text{for } r \leq r^* \\
        - \left(\frac{p}{p-1}\right) r^{*\frac{1}{p-1}} + (1 - r^*)^{(1/(p-1))-1} (r - r^*) \\
        + 1 - r^{*(p/(p-1))} - (1-r^*)^{(p/(p-1))} + 1 - E(r^*) & \text{for } r > r^* \\
    \end{cases}
    \label{eq:exact_solution}
\end{equation}

\noindent where \begin{align}
    F(r) &= r^{\frac{p}{p-1}}, \\
    G(r) &= (1-r)^{\frac{p}{p-1}}, \\
    E(r) &= \frac{p}{p-1}r.
\end{align}

\noindent \( u_{exact} \) is continuous and differentiable (i.e., \( C^1 \)). Moreover, the value of \( u_{exact} \) can be freely chosen for \( r > r^* \) without violating its \( C^1 \) smoothness.
The obstacle function, exact solution, and the neural network's approximation are shown for \( p=3 \) ( Figure \ref{fig:solution_plot_p_3}):

\begin{figure}[H]
    \centering
    \includegraphics[width=\linewidth]{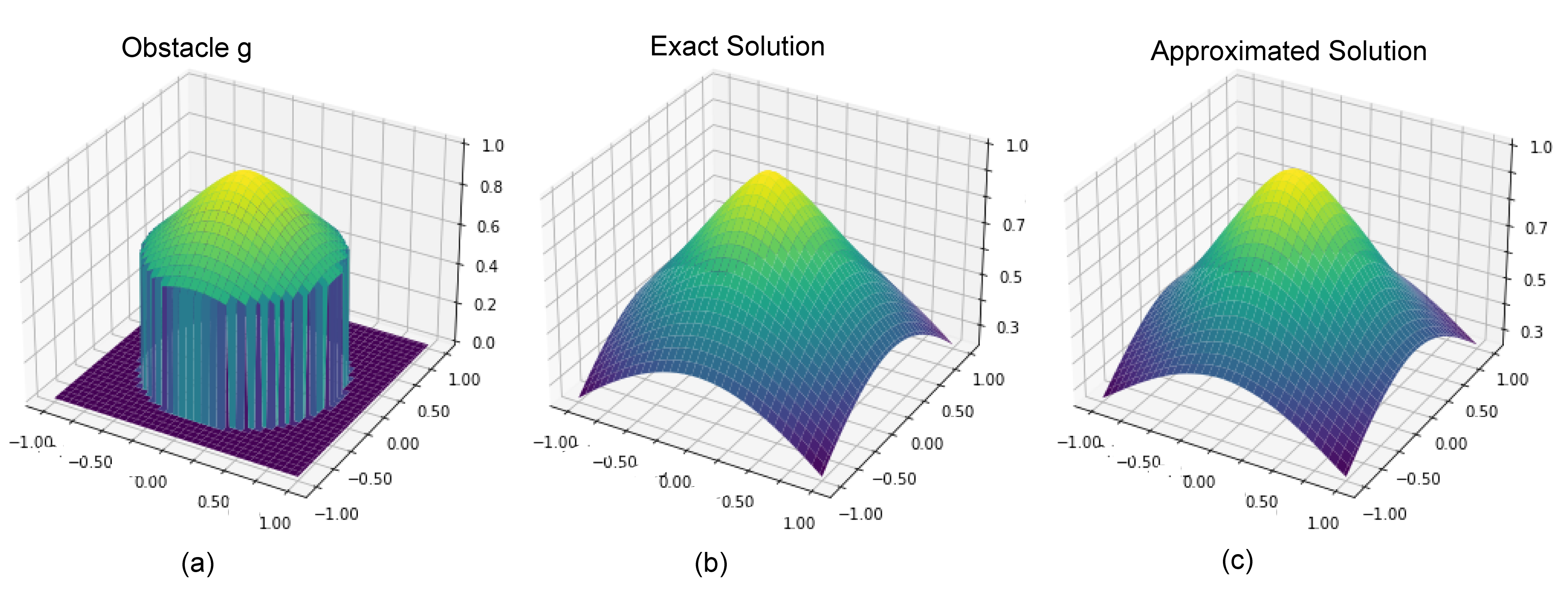}
    \caption{(a) Obstacle function plot, (b) Exact solution plot, and (c) Approximated solution plot.}
    \label{fig:solution_plot_p_3}
\end{figure}

\begin{figure}[H]
    \centering
    \includegraphics[width= \linewidth]{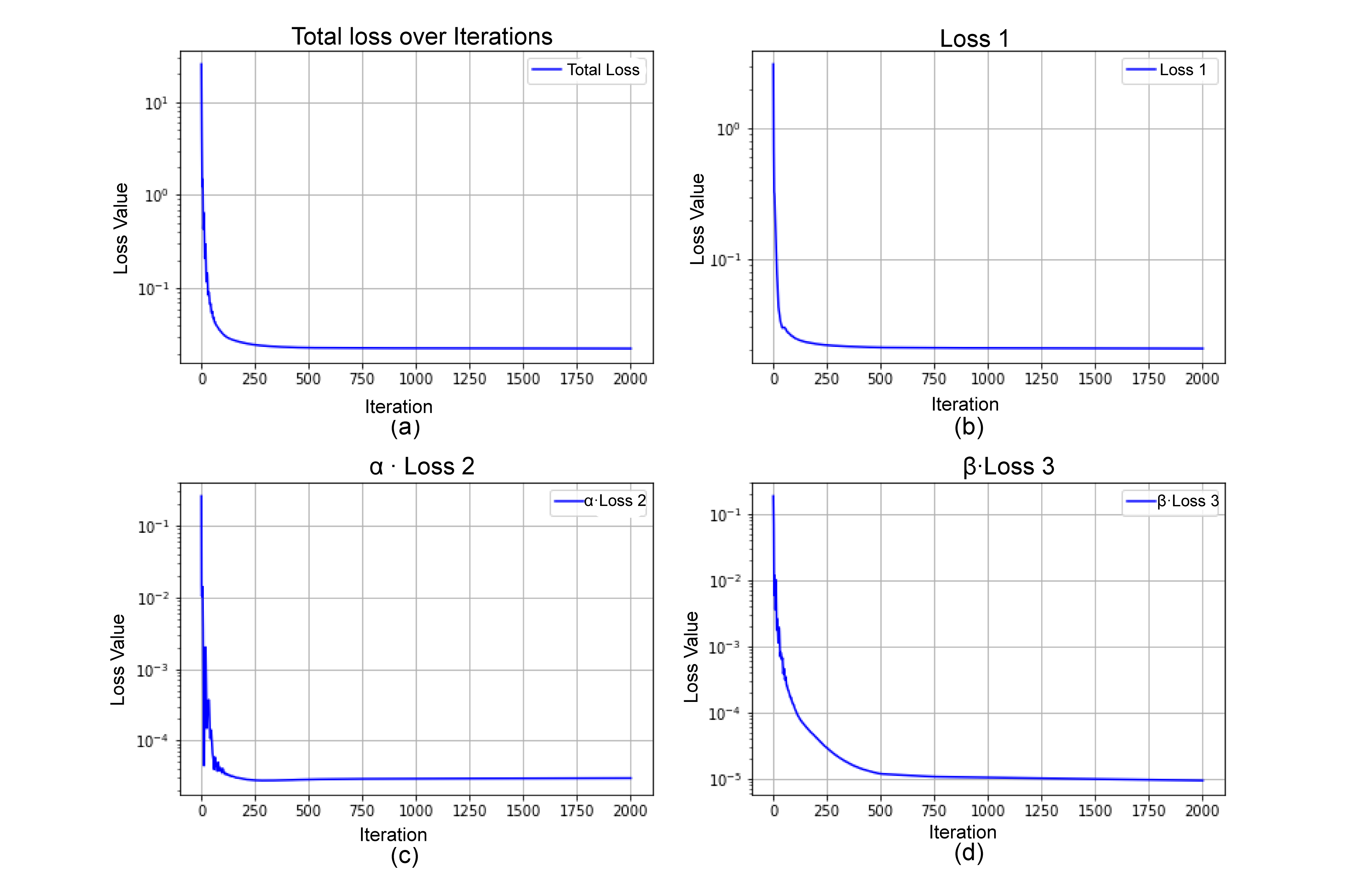}
    \caption{Evolution of different losses during training plotted on a log scale versus the number of iterations. (a) \textit{Total Loss}, as described in Equation \ref{total_loss}. (b) \textit{Loss 1} from Equation \ref{loss_1}. (c) \textit{Loss 2} influenced by coefficient $\alpha$, from Equation \ref{loss_2}. (d) \textit{Loss 3}, affected by coefficient $\beta$, as per Equation \ref{loss_3}.}
    \label{fig:lossplots_p_3}
\end{figure}

\begin{figure}[H]
    \centering
    \includegraphics[width=0.6\textwidth]{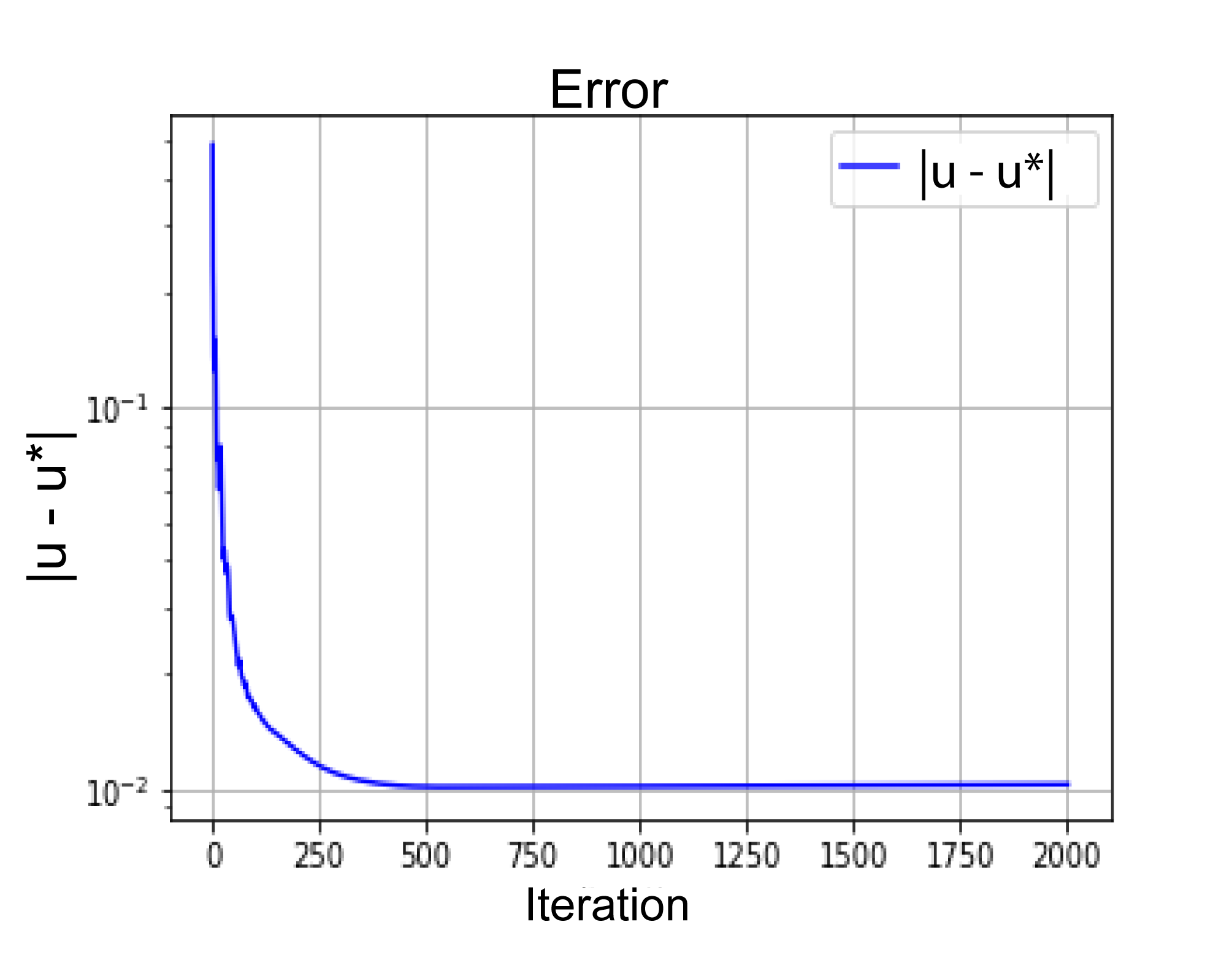}
    \caption{Evolution of the \(L^1\) error for \( p=3 \) between the approximated solution \( u^* \) and the exact solution \( u \) over training iterations. As with the previous case, this plot gives insight into the neural network's convergence behavior and prediction accuracy for this parameter setting.}
    \label{fig:exact_loss_plot_p_3}
\end{figure}

In this study, we set weighting parameters \(\alpha\) and \(\beta\) both to \(100\) and train the network for 2000 iterations. All losses are shown on a logarithmic scale against the number of iterations to illustrate the efficiency of the training process (\ref{fig:lossplots_p_3}). The $L^1$ norm of the difference between the approximated and true solutions as a function of training iteration is further shown in (\ref{fig:exact_loss_plot_p_3}). 

\begin{figure}[H]
    \centering
    \includegraphics[width=\linewidth]{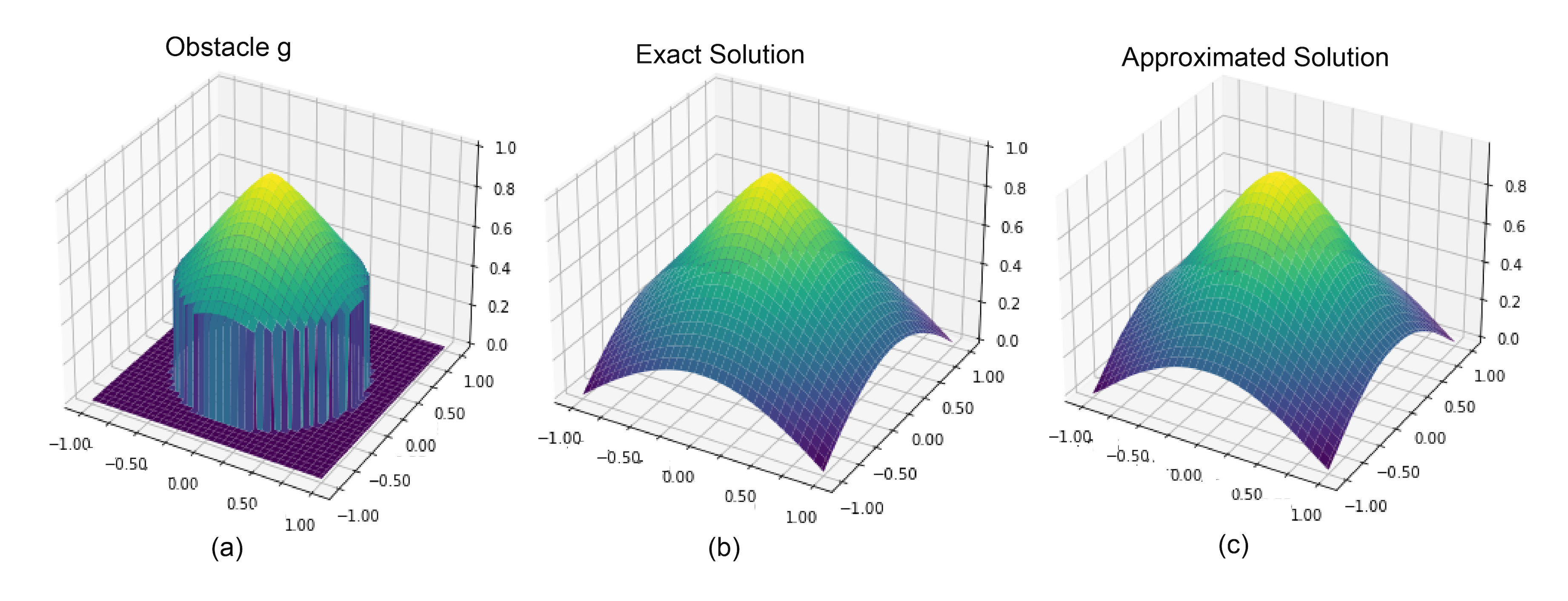}
    \caption{(a) Obstacle function plot, (b) Exact solution plot, and (c) Approximated solution plot.}
    \label{fig:solution_plot_p_4}
\end{figure}

For the case \( p = 4 \), Building upon our previous methodologies, we present the results for this scenario. Visual representations of the obstacle function, exact solution, and the neural network's approximation for \( p = 4 \) are shown in the figure (\ref{fig:solution_plot_p_4}).

\newpage

Post-training, we juxtapose the exact solution and the neural network's approximation for a more comprehensive visual comparison for \( p = 4 \) (Figure \ref{fig:solution_plot_p_4}). Similar to the previous cases, we showcase the loss behavior during training for \( p = 4 \) in the Figure \ref{fig:lossplots_p_4}. To assess the efficiency of the neural network's predictions for \( p = 4 \), we compute the \(L^1\) error (Figure \ref{fig:exact_loss_plot_p_4}):

\begin{figure}[H]
    \centering
    \includegraphics[width= \linewidth]{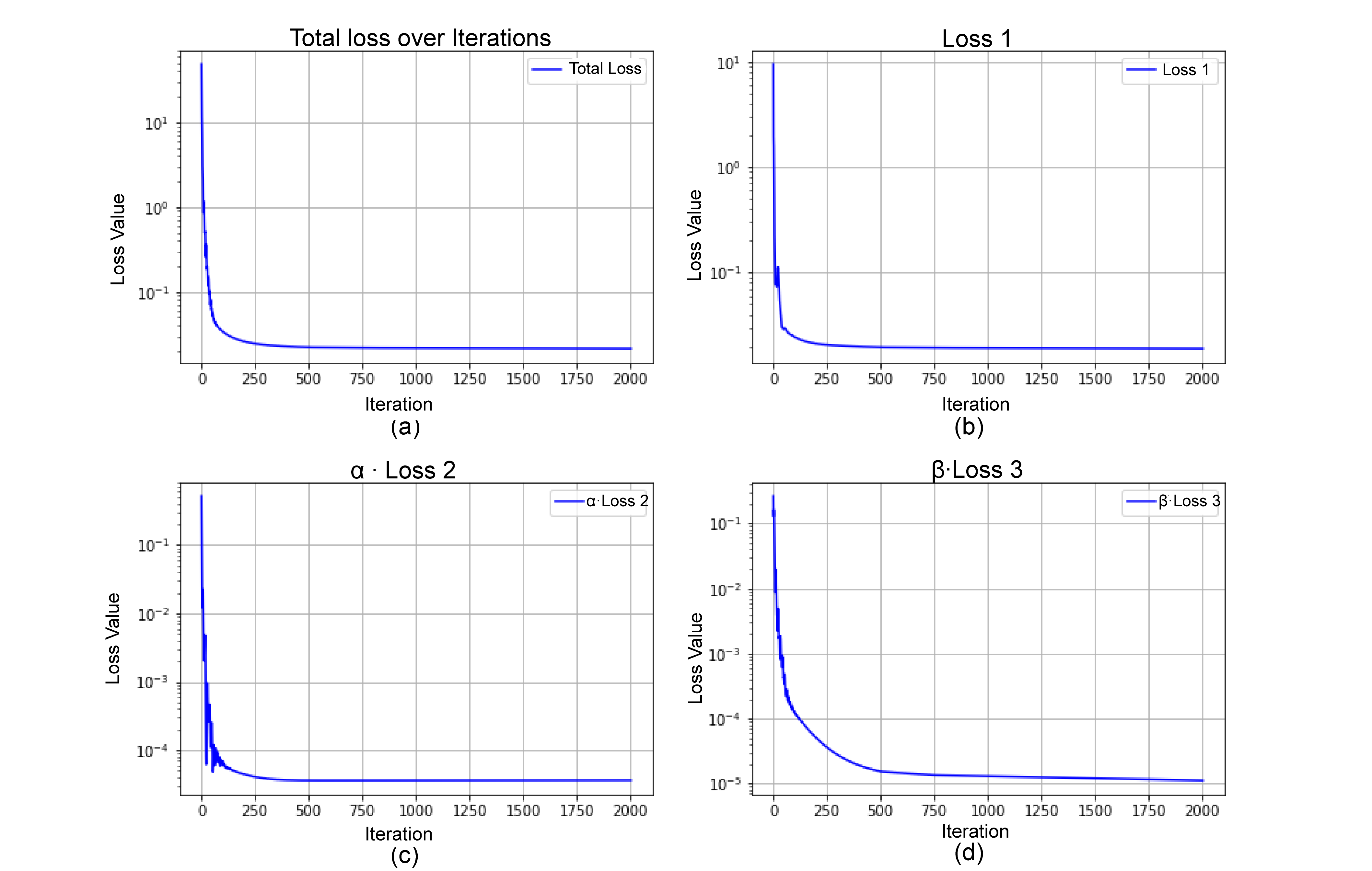}
    \caption{Evolution of different losses during training plotted on a log scale versus the number of iterations. (a) \textit{Total Loss}, as described in Equation \ref{total_loss}. (b) \textit{Loss 1} from Equation \ref{loss_1}. (c) \textit{Loss 2} influenced by coefficient $\alpha$, from Equation \ref{loss_2}. (d) \textit{Loss 3}, affected by coefficient $\beta$, as per Equation \ref{loss_3}.}
    \label{fig:lossplots_p_4}
\end{figure}

\begin{figure}[H]
    \centering
    \includegraphics[width=0.6\textwidth]{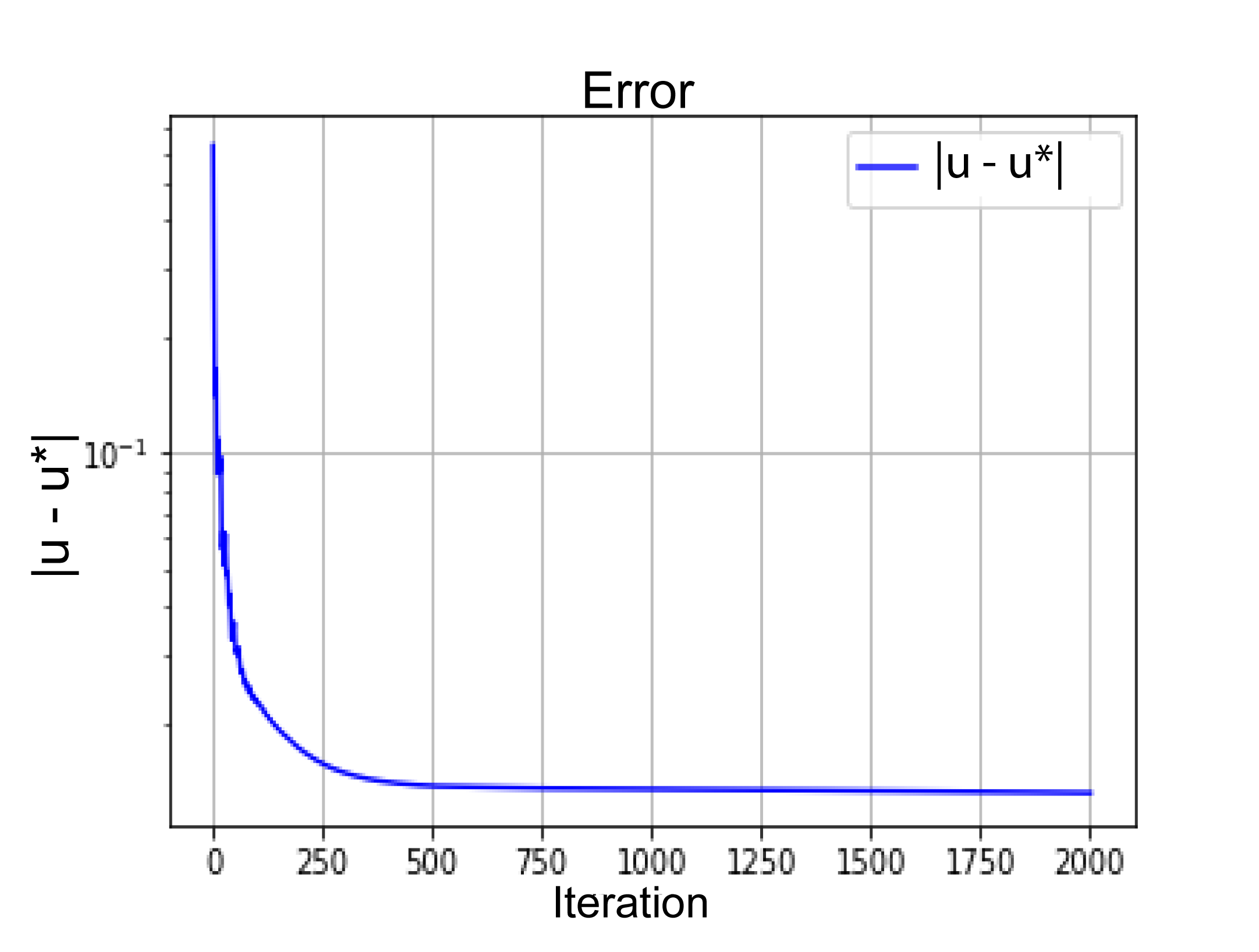}
    \caption{Evolution of the \(L^1\) error for \( p=4 \) between the approximated solution \( u^* \) and the exact solution \( u_{\text{exact}} \) over iterations, highlighting the convergence and prediction accuracy.}
    \label{fig:exact_loss_plot_p_4}
\end{figure}

These results demonstrate this approach can well approximate the solution of the obstacle constrained ice sheet PDE in simplified 1D and 2D domains. We next apply this approach to solve this PDE with a realistic obstacle geometry: the Greenland bedrock elevation taken from real-world data. The challenges and nuances of this dataset offer a rigorous testbed to evaluate the adaptability and robustness of our proposed algorithms.

\section{Application to Greenland Data}

In this section, we apply the techniques developed earlier to a dataset derived from Greenland. This dataset, available from the National Snow and Ice Data Center (NSIDC), offers a comprehensive view of various parameters related to Greenland's ice melt and movement. Specifically, we utilize data from the [NSIDC-0092 dataset \cite{Bamber2001}](\url{https://nsidc.org/data/nsidc-0092/versions/1}), which captures detailed information on ice-thickness. The NSIDC-0092 dataset offers detailed insights into Greenland's topography, presenting Digital Elevation Models (DEMs), ice thickness, and bedrock elevation data. The parameter range is extensive: DEM values range from -0.1 m to 3278.3 m, ice thickness measurements extend from 0 m to 3366.5 m, and bedrock elevation data vary between -963.1 m and 3239.0. These data are organized into ASCII text grids with a spatial resolution of 5 km, comprising 301 columns by 561 rows. This grid resolution ensures that the dataset provides a comprehensive and granular view of the ice sheet and bedrock topography across Greenland, suitable for precise climatological and glaciological analyses. The interpolation to a 5 km grid underscores the dataset's utility in modeling and analyzing the dynamics of Greenland's ice sheet and underlying bedrock with considerable detail.

To provide context, we begin by visualizing the bedrock topography from the dataset (\ref{fig:bedrock_topography_plot}):

\begin{figure}[H]
    \centering
    \includegraphics[width=0.7\textwidth]{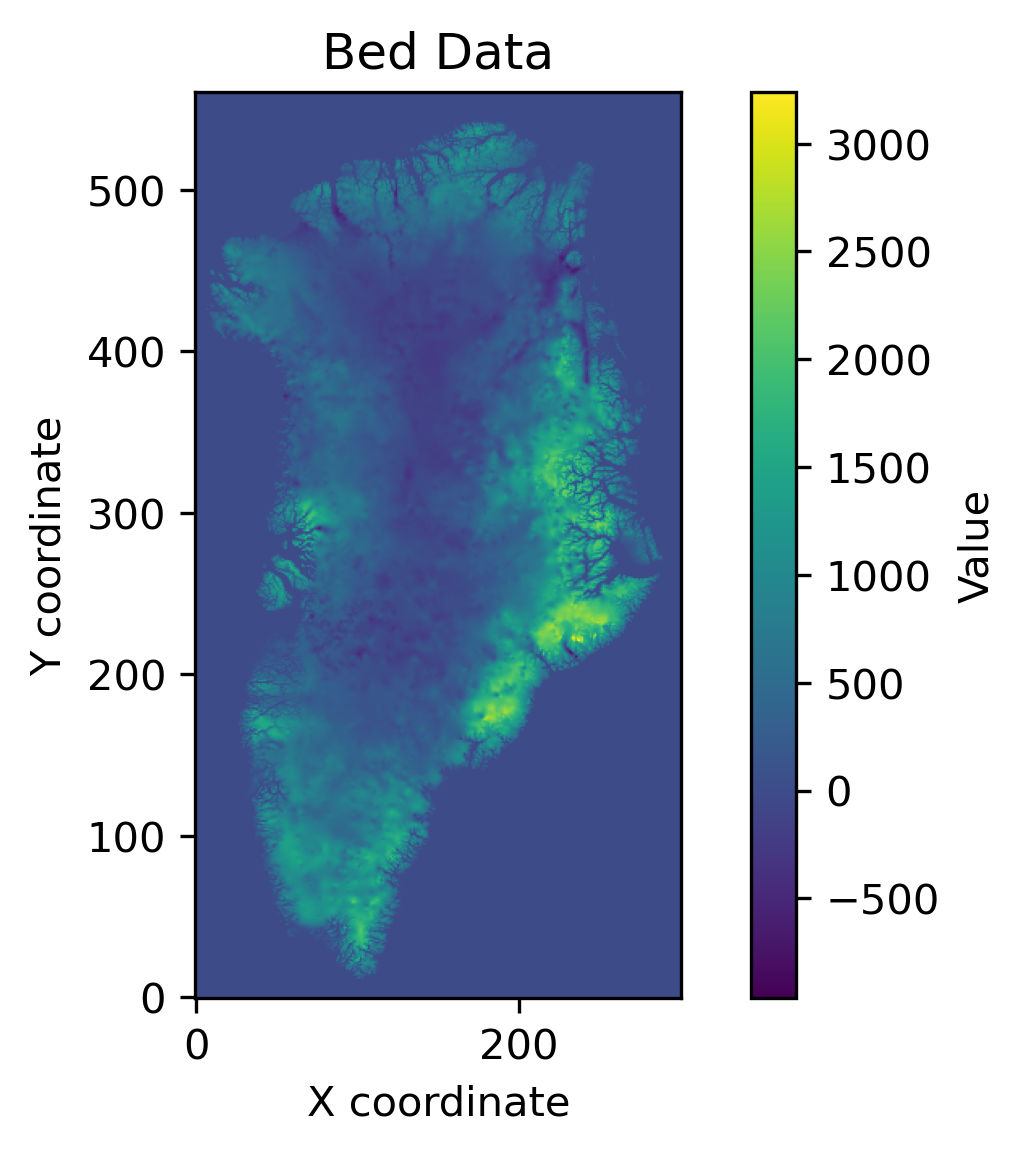}
    \caption{Bedrock topography of Greenland derived from the NSIDC-0092 dataset. This representation serves as our obstacle function for the subsequent analysis.}
    \label{fig:bedrock_topography_plot}
\end{figure}

\subsection{Neural Network Initialization and Training}
Translating this method to the bedrock obstacle data from the Greenland's ice-thickness data introduces a specific challenge. Notably, when employing conventional random initialization, the neural network behaved erratically, at least in part because it failed to satisfy the obstacle, leading to numerical issues. The resulting solutions were not only characterized by high error rates but were also noticeably distant from the anticipated results. Such deviations and unpredictable behavior underscored the need for an alternative approach to initialize the neural network.

To overcome this, we took a two-stage modeling approach. In the first, we intelligently initialize the NN. In the second, we train that initialized network to approximate the PDE (the same method as previously).

To initialize the network, we trained it to mimic the bedrock elevation as illustrated in Figure \ref{fig:bedrock_topography_plot}. This preliminary training phase enabled the neural network to learn and approximate the patterns of the bedrock topography. This is not a solution to the problem, only a method of initializing the network.

\begin{figure}[H]
    \centering
    \includegraphics[width=0.65\textwidth]{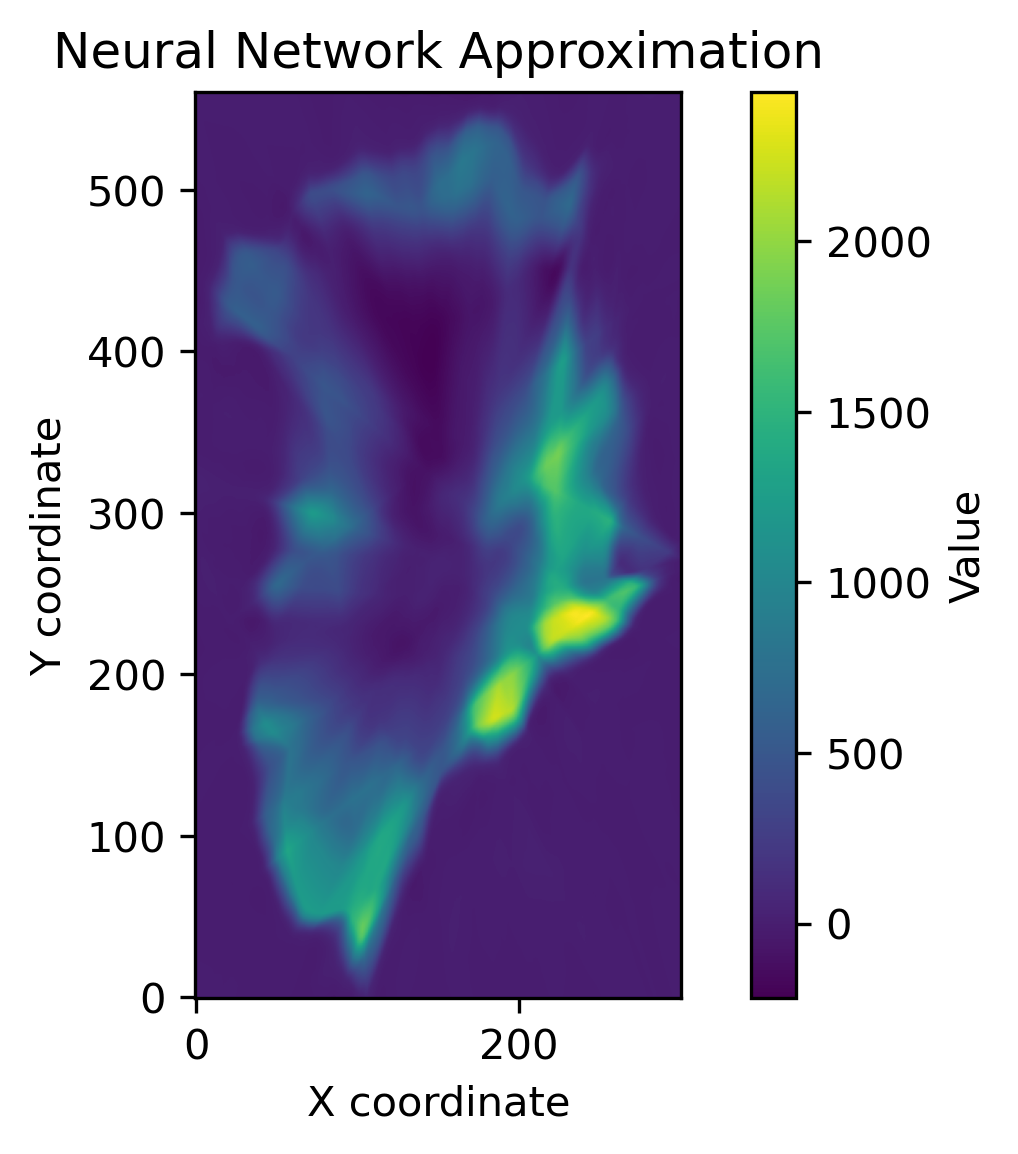}
    \caption{Visualization of the neural network after pre-training on the bedrock elevation. This pre-training produces an initialization for the network that is subsequently used as a starting point for training..}
    \label{fig:pretrained_network_visualization}
\end{figure}

Once this stage achieved satisfactory convergence, we retained the optimized weights of the network. These pre-trained weights were then employed as initialization parameters. By harnessing the pre-trained weights, the neural network exhibited a more efficient learning process and showed considerable improvements in error reduction. To summarize, we pre-trained the network to mimic a reasonable starting function, then used a second phase of training (same as previously discussed) to construct a solution to the PDE.

For the training process, we set the penalty parameters as \( \alpha = \beta = 4000 \). This was chosen by trail and error to ensure all weighted components of the loss are the same order of magnitude. For the solution approximation, we adopt an architecture comprising 15 hidden layers, each equipped with 320 neurons. The neural network underwent training for a total of \( 22000 \) iterations to ensure accurate convergence and approximation of the Greenland ice-thickness data ( \ref{fig:loss_greenland_p3}).

\begin{figure}[H]
    \centering
    \includegraphics[width=\linewidth]{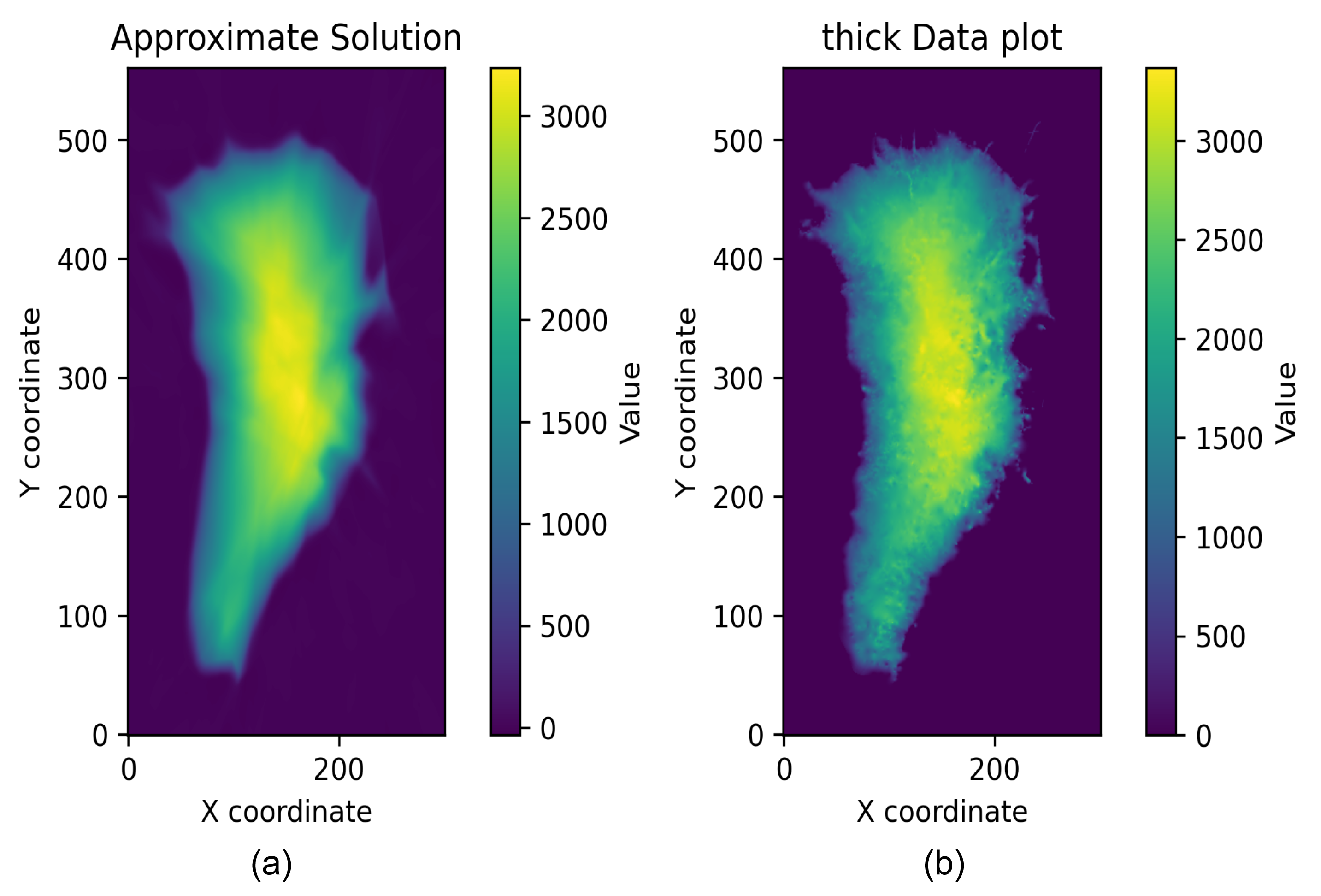}
    \caption{Comparison over Greenland: (a) Neural network-derived approximation for \(p=4\) and (b) precise data from the NSIDC-0092 dataset.}
    \label{fig:greenland_comparison}
\end{figure}

Following the training process and the application of our methodology, we obtained the approximated solution. For a clear comparison, we juxtapose this solution alongside the exact ice thickness data derived from the NSIDC-0092 dataset (\ref{fig:greenland_comparison}).

\subsection{Analysis of Training Losses}

\begin{figure}[H]
    \centering
    \includegraphics[width= \linewidth]{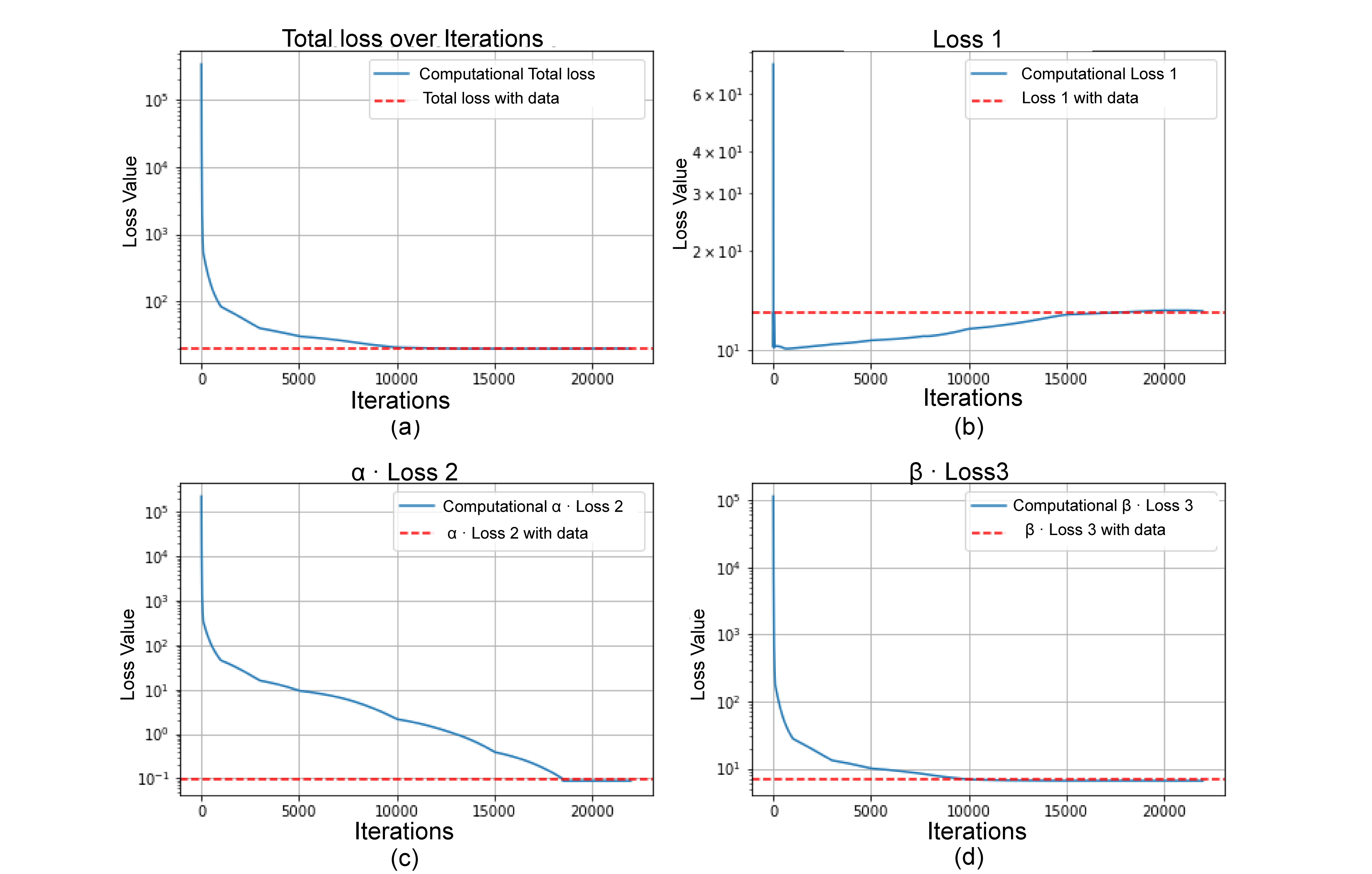}
    \caption{Comparison of various loss components for \( p=3 \). Each plot showcases the loss trajectory over iterations, with the red dashed line representing the loss upon data input, and the blue solid line symbolizing computational loss.}
    \label{fig:loss_greenland_p3}
\end{figure}

\begin{figure}[H]
    \centering
    \includegraphics[width= \linewidth]{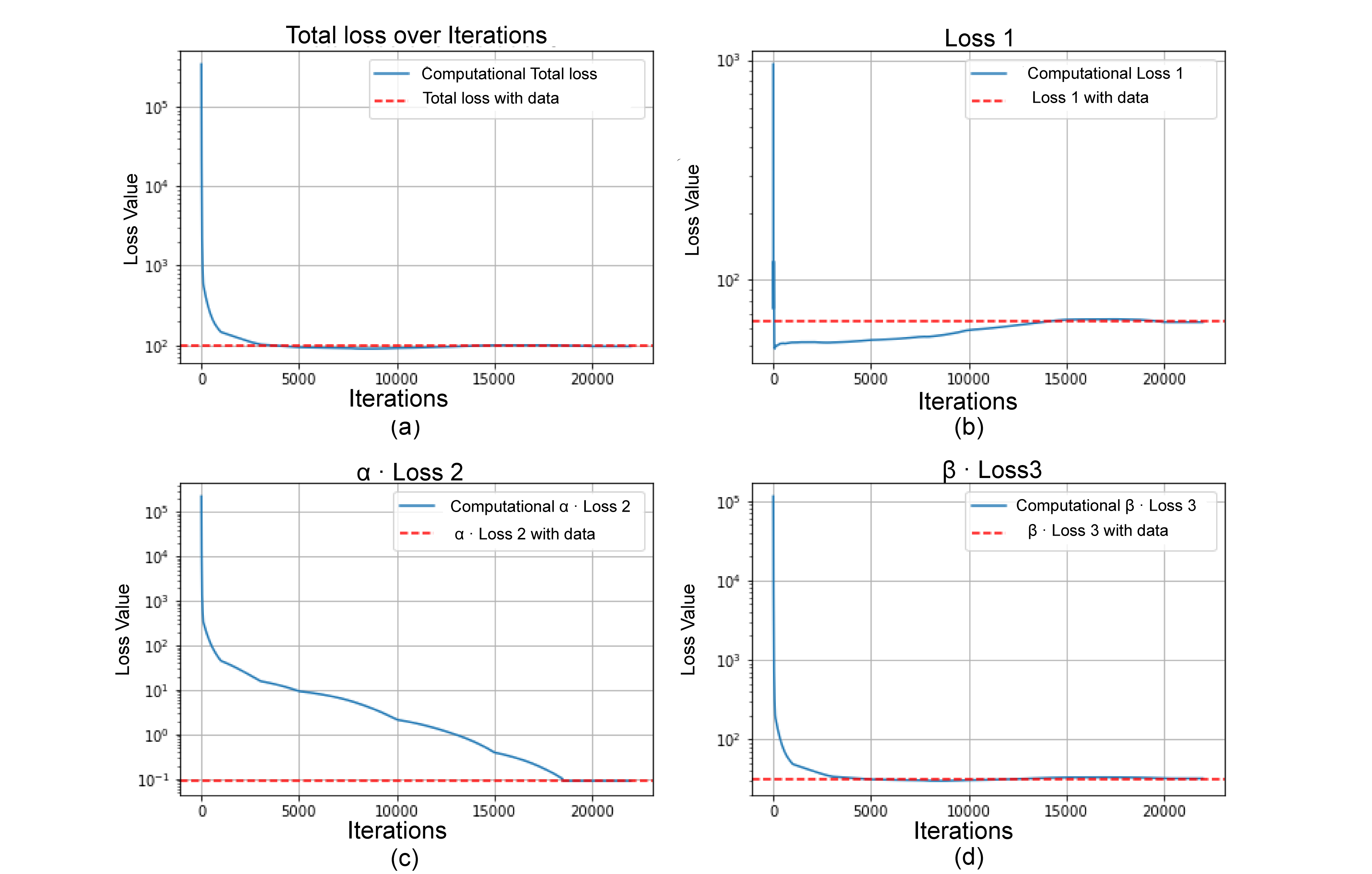}
    \caption{Comparison of various loss components for \( p=4 \). Each plot showcases the loss trajectory over iterations, with the red dashed line representing the loss upon data input, and the blue solid line symbolizing computational loss.}
    \label{fig:loss_greenland_p4}
\end{figure}

Determining how well the approximated solution does at satisfying the PDE is challenging for two reasons: 1) there is no exact solution to compare against and 2) the magnitude of training losses are not interpretable in absolute terms.  To circumvent this issue, we compare the properties of the neural network solution to this problem with the ice-sheet thickness data. Note however that the thickness data is itself not necessarily a solution to the PDE, a research question that is beyond the scope of this article. Thus we will not directly compare the solution to the data. Instead, to get a benchmark for comparison, we compute the losses that would be found from inserting the measured ice sheet data into the three loss terms. These form the dashed lines against which training losses are compared in Figures (\ref{fig:loss_greenland_p3}) and  (\ref{fig:loss_greenland_p4}). These dashed lines provide a scale reference to compare the model training losses against.

Figure (\ref{fig:loss_greenland_p3}) shows the trajectory of the losses as a function of training iteration. Results show that all losses converge to an absolute level that is consistent with the real ice sheet thickness data's satisfaction of the PDE. This, in combination with the visual comparison of the approximated solution to the real data suggests this approach well approximates the solution to this PDE. Similar analysis was carried out for $p=4$ with similar results (Figures (\ref{fig:loss_greenland_p4}))

\subsection{Loss Analysis for Different \texorpdfstring{\( p \) Values}{p Values}}

\section{Summary}

This paper embarked on an exploration of the integration between the realm of mathematical modeling and deep learning to address intricate problems related to the dynamics of Greenland's ice sheet.
 We first established a robust theoretical framework, illustrating the underpinnings of the variational inequalities and obstacle problems. Adopting a novel approach, we leveraged the power of neural networks to seek solutions, demonstrating both its theoretical and practical applications.

A key highlight of our methodology was the integration of pre-training. By first acquainting the neural network with a representation of the bedrock topography, we ensured a more efficient learning process when the actual ice-thickness data was introduced. This strategic move not only optimized our results but also underscores the importance of judicious model initialization in complex problem-solving.

When applied to the NSIDC-0092 dataset, our method displayed remarkable accuracy. The approximated solutions closely mirrored the exact data, which stands as a testament to the promise and potential of this interdisciplinary approach. 

Furthermore, this study illuminates the broader implications for climate science. As the melting of Greenland's ice sheet plays a pivotal role in global climate dynamics, having precise and efficient computational methods becomes paramount. The methodologies explored here can be adapted and extended to other similar domains, making it a versatile tool in the scientific community's arsenal.

In this paper, we've shown how deep learning can work together with mathematical models to study Greenland's ice sheet dynamics. As we collect more data and face more complex challenges in the future, using both these tools together will help us find better solutions. This combination is promising for future research and problem-solving.

\section{Acknowledgements}
The authors acknowledge useful discussions on the mathematical background with Professor Ed. Bueler and Paolo Piersanti. This work was supported in part by the Research Fund of Indiana University.

% References
\bibliographystyle{plain}
\bibliography{paper}

\end{document}